\documentclass{article}
\usepackage{amsfonts}

\input{tcilatex}
\begin{document}

\title{\bigskip \bigskip BOUNDS FOR THE DISTANCE ESTRADA INDEX OF GRAPHS}
\author{\c{S}. Burcu Bozkurt and Durmu\c{s} Bozkurt \\
\textit{Department of Mathematics, Science Faculty, }\\
\textit{Sel\c{c}uk University, 42075, Campus, Konya, Turkey}\\
\textit{sbbozkurt@selcuk.edu.tr, dbozkurt@selcuk.edu.tr}}
\maketitle

\begin{abstract}
The $D$-eigenvalues $\mu _{1},\mu _{2},\ldots ,\mu _{n}$ of a connected
graph $G$ are the \linebreak eigenvalues of its distance matrix. The
distance Estrada index of $G$ is defined in [15] as 
\[
DEE=DEE\left( G\right) =\dsum\limits_{i=1}^{n}e^{\mu _{i}} 
\]

In this paper, we give better lower bounds for the distance Estrada index of
any connected graph as well as some relations between $DEE\left( G\right) $
and the distance energy.


\textit{2000 AMS Classification: }05C12, 05C90.
\end{abstract}

\bigskip

\section{Introduction}

Let $G$ be a simple connected graph with $n$ vertices and $m$ edges on
vertex set $V\left( G\right) =\left\{ v_{1},v_{2},\ldots ,v_{n}\right\} .$
Throughout this paper, such a graph will be referred to as connected $\left(
n,m\right) $-graph. The distance matrix $D=D\left( G\right) $ of $G$ is
defined so that its $\left( i,j\right) $-entry is equal to $d_{G}\left(
v_{i},v_{j}\right) $, denoted by $d_{ij}$, the distance (i.e., the length of
the shortest path [1]) between the vertices $v_{i}$ and $v_{j}$ of $G$. The
diameter of the graph $G$ is the maximum distance between any two vertices
of $G$. Let $\Delta $ be diameter of $G$ and $A\left( G\right) $ be the
adjacency matrix of $G$. The eigenvalues of $A\left( G\right) $ are called
the eigenvalues of $G$ and the eigenvalues of $D\left( G\right) $ are said
to be the $D$-eigenvalues of $G$ [2]. Since $A\left( G\right) $ and $D\left(
G\right) $ are real symmetric matrices their eigenvalues are real numbers.
So we can order them so that $\lambda _{1}\geq \lambda _{2}\geq \cdots \geq
\lambda _{n}$ and $\mu _{1}\geq \mu _{2}\geq \cdots \geq \mu _{n}$ are the
eigenvalues and $D$-eigenvalues of $G$, respectively.

\bigskip

The Estrada index of the graph $G$ is defined in [4-8] as:%
\begin{equation}
EE=EE\left( G\right) =\dsum\limits_{i=1}^{n}e^{\lambda _{i}}  \tag{1}
\end{equation}

The Estrada index of graphs has an important role in Chemistry and Physics.
There exist a vast literature related to Estrada index and its bounds. For
more information see [4-12]. Because of the evident success of the graph
Estrada index, Estrada index based of the eigenvalues of other graph
matrices have, one-by-one, been introduced. From this respect, the authors
defined the distance Estrada index $DEE\left( G\right) $ based on distance
matrix of the graph $G$ as the following

\begin{equation}
DEE=DEE\left( G\right) =\dsum\limits_{i=1}^{n}e^{\mu _{i}}.  \tag{2}
\end{equation}%
where $\mu _{1},\mu _{2},\ldots ,\mu _{n}$ are the $D$-eigenvalues of $G.$

Let 
\[
N_{k}=\dsum\limits_{i=1}^{n}\left( \mu _{i}\right) ^{k}. 
\]

Recalling the power series expansion of $e^{x}$ we have another expression
of distance Estrada index as the following

\begin{equation}
DEE\left( G\right) =\dsum\limits_{k=0}^{\infty }\frac{N_{k}}{k!}.  \tag{3}
\end{equation}

Some well known mathematical properties on the distance Estrada index of the
graph $G$ are established as follows:

\bigskip

\textbf{\ Theorem 1.1}\textit{\ [15]\ Let }$G$\textit{\ be a connected }$%
\left( n,m\right) $\textit{-graph and }$\Delta $\textit{\ the diameter of }$%
G $\textit{. Then the distance Estrada index is bounded as follows}%
\begin{equation}
\sqrt{n^{2}+4m}\leq DEE\left( G\right) \leq n-1+e^{\Delta \sqrt{n\left(
n-1\right) }}.  \tag{4}
\end{equation}%
\textit{Equality holds on both sides of (4) if and only if }$G\simeq K_{1}.$

\bigskip

\textbf{Theorem 1.2 }\textit{[15] Let }$G$\textit{\ be a connected }$\left(
n,m\right) $\textit{-graph and }$\Delta $\textit{\ the diameter of }$G$%
\textit{. Then}%
\begin{equation}
DEE\left( G\right) -E_{D}\left( G\right) \leq n-1-\Delta \sqrt{n\left(
n-1\right) }+e^{\Delta \sqrt{n\left( n-1\right) }}  \tag{5}
\end{equation}%
\textit{or}%
\begin{equation}
DEE\left( G\right) \leq n-1+e^{E_{D}\left( G\right) }.  \tag{6}
\end{equation}%
\textit{Equality holds in (5) or (6) if and only if }$G\simeq K_{1}.$

\bigskip

We organize this paper in the following way. In section 2, we obtain
\linebreak better lower bounds for $DEE\left( G\right) $. In section 3, we
give some relations between $DEE\left( G\right) $ and the distance energy.

\bigskip

\section{\textbf{Lower bounds for the distance Estrada index}}

\bigskip

In this section, we obtain some lower bounds for the distance Estrada index
of any connected graph $G$ using the same procedure in [12]. Firstly, we
give the following definitions and lemmas which will be needed then.

\bigskip \textbf{\ Definition 2.1 }[14] Let $G$ be a graph with $V\left(
G\right) =\left\{ v_{1},v_{2},\ldots ,v_{n}\right\} $ and a distance matrix $%
D$. Then the distance degree of $v_{i}$, denoted by $D_{i}$ is given by $%
D_{i}=\dsum\limits_{j=1}^{n}d_{ij}.$

\bigskip \textbf{Definition 2.2 }[14]\textbf{\ }Let $G$ be a graph with $%
V\left( G\right) =\left\{ v_{1},v_{2},\ldots ,v_{n}\right\} $ and a distance
matrix $D$. Let the distance degree sequence be $\left\{ D_{1},D_{2},\ldots
,D_{n}\right\} .$ Then the second distance degree of $v_{i}$, denoted by $%
T_{i}$ is given by $T_{i}=\dsum\limits_{j=1}^{n}d_{ij}D_{j}.$

\bigskip \textbf{Definition 2.3 } [16]\textbf{\ }Let $G$ be a graph with $%
V\left( G\right) =\left\{ v_{1},v_{2},\ldots ,v_{n}\right\} $ and a distance
matrix $D$. Let the distance degree sequence be $\left\{ D_{1},D_{2},\ldots
,D_{n}\right\} .$Then for each $i=1,2,\ldots ,n$ the sequence $%
M_{i}^{(1)},M_{i}^{(2)},\ldots ,M_{i}^{(t)},\ldots $ \ is defined as
follows: Fix $\alpha \in 
\mathbb{R}
$, let 
\[
M_{i}^{(1)}=D_{i}^{\alpha } 
\]%
and for each $t\geq 2$, let 
\[
M_{i}^{(t)}=\dsum\limits_{j=1}^{n}d_{ij}M_{j}^{(t-1)}. 
\]

\textbf{Definition 2.4} [3] Let $G$ be a graph with distance matrix $D$.
Then the Wiener index of $G$, denoted by $W\left( G\right) $ is given by $%
W\left( G\right) =\frac{1}{2}\dsum\limits_{i=1}^{n}D_{i}.$\textbf{\ }

\textbf{Lemma 2.1} \textit{[16] Let }$G$\textit{\ be a connected graph }$%
\alpha $\textit{\ be a real number and }$t$\textit{\ be an integer }$.$%
\textit{Then}%
\[
\mu _{1}\geq \sqrt{\frac{S_{t+1}}{S_{t}}}. 
\]%
\textit{where }$S_{t}\mathit{=}\dsum\limits_{i=1}^{n}\left(
M_{i}^{(t)}\right) ^{2}.$ \textit{Moreover, equality holds for particular
values of }$\alpha $\textit{\ and }$t$\textit{\ if and only if }%
\[
\frac{M_{1}^{(t+1)}}{M_{1}^{(t)}}=\frac{M_{2}^{(t+1)}}{M_{2}^{(t)}}=\cdots =%
\frac{M_{n}^{(t+1)}}{M_{n}^{(t)}}. 
\]

\textbf{Lemma 2.2}\textit{\ [14]\ A connected graph }$G$\textit{\ has two
distinct }$D$\textit{-eigenvalues if and only if }$G$\textit{\ is a complete
graph.}

\bigskip

\textbf{Lemma 2.3}\textit{\ [17]\ Let }$G$\textit{\ be a connected graph
with }$n\geq 2$\textit{\ vertices and }$m$\textit{\ edges. Then}%
\[
\mu _{1}\geq 2\left( n-1\right) -\frac{2m}{n} 
\]%
\textit{with equality if and only if }$G=K_{n}$\textit{\ or }$G$\textit{\ is
a regular graph of diameter two.}

\bigskip

Now we are ready to give the main results of this section.

\bigskip

\textbf{Theorem 2.1 }\textit{Let }$G$\textit{\ be a connected }$\left(
n,m\right) $-\textit{graph. Then}%
\begin{equation}
DEE(G)\geq e^{\sqrt{\frac{S_{t+1}}{S_{t}}}}+\frac{n-1}{e^{\frac{1}{n-1}\sqrt{%
\frac{S_{t+1}}{S_{t}}}}}.  \tag{7}
\end{equation}%
\textit{where }$\alpha $\textit{\ is a real number, }$t$\textit{\ is an
integer and }$S_{t}\mathit{=}\dsum\limits_{i=1}^{n}\left( M_{i}^{(t)}\right)
^{2}.$\textit{\ Moreover, the equality holds in (7) if and only if }$G$%
\textit{\ is the complete graph }$K_{n}.$

\textit{\bigskip }

\textit{Proof.} Starting with the equation (2) and using
Arithmetic-Geometric Mean Inequality, we get%
\begin{eqnarray}
DEE(G) &=&e^{\mu _{1}}+e^{\mu _{2}}+\cdots +e^{\mu _{n}}  \nonumber \\
&\geq &e^{\mu _{1}}+(n-1)\left( \dprod\limits_{i=2}^{n}e^{\mu _{i}}\right) ^{%
\frac{1}{n-1}}  \TCItag{8} \\
&=&e^{\mu _{1}}+(n-1)\left( e^{-\mu _{1}}\right) ^{\frac{1}{n-1}},\text{
since }\dsum\limits_{i=1}^{n}\mu _{i}=0.  \TCItag{9}
\end{eqnarray}%
Consider the following function 
\[
f(x)=e^{x}+\frac{n-1}{e^{\frac{x}{n-1}}} 
\]%
for $x>0$. We have 
\[
f%
{\acute{}}%
(x)=e^{x}-e^{-\frac{x}{n-1}}>0 
\]%
for $x>0$. It is easy to see that $f$ is an increasing function for $x>0$.
From the equation (9) and Lemma 2.1, we obtain

\begin{equation}
DEE(G)\geq e^{\sqrt{\frac{S_{t+1}}{S_{t}}}}+\frac{n-1}{e^{\frac{1}{n-1}\sqrt{%
\frac{S_{t+1}}{S_{t}}}}}.  \tag{10}
\end{equation}%
This completes the first part of the proof.

\bigskip

Now we suppose that the equality holds in (7). Then all inequalities in the
above argument must be equalities. From (10) we have%
\[
\mu _{1}=\sqrt{\frac{S_{t+1}}{S_{t}}} 
\]%
which implies $\frac{M_{1}^{(t+1)}}{M_{1}^{(t)}}=\frac{M_{2}^{(t+1)}}{%
M_{2}^{(t)}}=\cdots =\frac{M_{n}^{(t+1)}}{M_{n}^{(t)}}$ $.$ From (8) and
Arithmetic-Geometric Mean Inequality we get $\mu _{2}=\mu _{3}=\cdots =\mu
_{n}.$ Therefore $G$ has exactly two distinct $D$-eigenvalues, by Lemma 2.2, 
$G$ is the complete graph $K_{n}.$

\bigskip Conversely, one can easily see that the equality holds in (7) for
the complete graph $K_{n}.$ This completes the proof.

\bigskip

The following result states a lower bound for the distance Estrada index
involving Wiener index.

\bigskip

\textbf{Corollary 2.1 }\textit{Let }$G$\textit{\ be a} connected $\left(
n,m\right) $-graph\textit{. Then}%
\begin{equation}
DEE(G)\geq e^{\frac{2W\left( G\right) }{n}}+\frac{n-1}{e^{\frac{2W\left(
G\right) }{n(n-1)}}}  \tag{11}
\end{equation}%
\textit{where }$W\left( G\right) $\textit{\ denotes the Wiener index of the
graph }$G$\textit{. Moreover the equality holds in (11) if and only if }$G$%
\textit{\ is the complete graph }$K_{n}.$

\textit{\bigskip }

\textit{Proof. }In [16], G\"{u}ng\"{o}r and Bozkurt showed that the
folllowing inequality (see Theorem 2)

\begin{equation}
\mu _{1}\geq \sqrt{\frac{S_{t+1}}{S_{t}}}\geq \sqrt{\frac{%
\dsum\limits_{i=1}^{n}T_{i}^{2}}{\dsum\limits_{i=1}^{n}D_{i}^{2}}}.  \tag{12}
\end{equation}%
where $S_{t}\mathit{=}\dsum\limits_{i=1}^{n}\left( M_{i}^{(t)}\right) ^{2}.$
Also in [14], Indulal proved the following result using Cauchy-Schwartz
inequality (see Theorem 4) 
\begin{equation}
\mu _{1}\geq \sqrt{\frac{\dsum\limits_{i=1}^{n}T_{i}^{2}}{%
\dsum\limits_{i=1}^{n}D_{i}^{2}}}\geq \sqrt{\frac{\dsum%
\limits_{i=1}^{n}D_{i}^{2}}{n}}\geq \frac{2W\left( G\right) }{n}.  \tag{13}
\end{equation}%
Combining Theorem 2.1, (12) and (13) we get the inequality (11).

Also, the equality holds in (11) if and only if $G$ is the complete graph $%
K_{n}.$

\bigskip

\textbf{Theorem 2.2 }\textit{\ Let }$G$\textit{\ be a connected }$\left(
n,m\right) $\textit{-graph with }$n\geq 2$\textit{\ vertices. Then}%
\begin{equation}
DEE(G)\geq e^{2\left( n-1\right) -\frac{2m}{n}}+e^{-\left( 2\left(
n-1\right) -\frac{2m}{n}\right) }+n-2.  \tag{14}
\end{equation}

\textit{Moreover the equality holds in (14) if and only if }$G=K_{2}$\textit{%
.}

\textit{\bigskip }

\textit{Proof. }Since $G$ is a connected graph $\mu _{1}>0$ and $\mu _{n}<0$%
. Now%
\begin{eqnarray}
DEE(G) &=&e^{\mu _{1}}+e^{\mu _{2}}+\cdots +e^{\mu _{n}}  \nonumber \\
&\geq &e^{\mu _{1}}+e^{\mu _{n}}+(n-2)\left( \dprod\limits_{i=2}^{n-1}e^{\mu
_{i}}\right) ^{\frac{1}{n-2}}  \TCItag{15} \\
&=&e^{\mu _{1}}+e^{\mu _{n}}+(n-2)e^{-\frac{\mu _{1}+\mu _{n}}{n-2}},\text{
since }\dsum\limits_{i=1}^{n}\mu _{i}=0.  \TCItag{16}
\end{eqnarray}%
We consider the following function%
\[
f(x,y)=e^{x}+e^{y}+(n-2)e^{-\frac{x+y}{n-2}} 
\]%
for $x>0$, $y<0$. In [12], Das and Lee showed that $f(x,y)$ has a minimum
value at $\ x+y=0$ and its minimum value is $e^{x}+e^{-x}+n-2$ (see Theorem
2.4)$.$ Also we can easily see that $e^{x}+e^{-x}+n-2$ is an increasing
function for $x>0.$ By Lemma 2.3 we obtain%
\begin{equation}
e^{\mu _{1}}+e^{-\mu _{1}}+n-2\geq e^{2\left( n-1\right) -\frac{2m}{n}%
}+e^{-\left( 2\left( n-1\right) -\frac{2m}{n}\right) }+n-2.  \tag{17}
\end{equation}%
From the equation (16) and (17), we get%
\[
DEE(G)\geq e^{2\left( n-1\right) -\frac{2m}{n}}+e^{-\left( 2\left(
n-1\right) -\frac{2m}{n}\right) }+n-2. 
\]%
This completes the first part of the theorem.

Now we suppose that the equality holds in (14). Then all inequalities in the
above argument must be equalities. Since $\mu _{1}+\mu _{n}=0$, we have that
-$\mu _{1}$ is also an $D$-eigenvalue of $G$. From equality in (17) and
Lemma 2.3 we can write $G$ is the complete graph $K_{n}$ $.$ From equality
in (15) and $\dsum\limits_{i=1}^{n}\mu _{i}=0$ we obtain%
\[
\mu _{2}=\mu _{3}=\cdots =\mu _{n-1}=0\text{ \ } 
\]%
since $\mu _{1}+\mu _{n}=0$. In [17], Zhou and Trinajsti\'{c} showed that
this situation is not valid for $n\geq 3$ (see Theorem 4)$.$ Therefore the
equality holds in (14) for only the graph $G=K_{2}.$

\bigskip Conversely, one can easily see that the equality holds in (14) for
\ the complete graph $K_{2}.$ This completes the proof.

\section{\protect\bigskip\ \textbf{Bounds for the distance Estrada index }%
\protect\linebreak \textbf{involving the distance energy}}

\bigskip

In this section, we firstly recall the distance energy $E_{D}\left( G\right) 
$ which is defined in [13] as 
\begin{equation}
E_{D}\left( G\right) =\dsum\limits_{i=1}^{n}\left\vert \mu _{i}\right\vert 
\tag{18}
\end{equation}%
where $\mu _{1},\mu _{2},\ldots ,\mu _{n}$ are the $D$-eigenvalues of $G$.

\bigskip

Now, we will adapt the some results in [10] on distance Estrada index to
give some relations between the distance Estrada index $DEE\left( G\right) $
and the distance energy $E_{D}\left( G\right) $.

\bigskip

\textbf{Theorem 3.1} \textit{Let }$G$\textit{\ be a connected }$\left(
n,m\right) $\textit{-graph.Then the distance Estrada index }$DEE\left(
G\right) $\textit{\ and the distance energy }$E_{D}\left( G\right) $\textit{%
\ satisfy the following inequality}%
\begin{equation}
\frac{1}{2}E_{D}\left( G\right) \left( e-1\right) +n-n_{+}\leq DEE\left(
G\right) \leq n-1+e^{\frac{E_{D}\left( G\right) }{2}}.  \tag{19}
\end{equation}%
\textit{where }$n_{+}$\textit{\ denotes the number of positive }$D$\textit{%
-eigenvalues of }$G.$\textit{\ Moreover, the equality holds on both sides of
(19) if and only if }$G=K_{1}.$

\bigskip \textit{Proof.}\textbf{\ Lower bound:}\textit{\ }Since $e^{x}\geq
ex,$ equality holds if and only if $x=1$ and $e^{x}\geq 1+x,$ equality holds
if and only if $x=0.$ We get

\begin{eqnarray*}
DEE\left( G\right) &=&\dsum\limits_{i=1}^{n}e^{\mu _{i}}=\dsum\limits_{\mu
_{i}>0}e^{\mu _{i}}+\dsum\limits_{\mu _{i}\leq 0}e^{\mu _{i}} \\
&\geq &\dsum\limits_{\mu _{i}>0}e\mu _{i}+\dsum\limits_{\mu _{i}\leq
0}\left( 1+\mu _{i}\right) \\
&=&e\left( \mu _{1}+\mu _{2}+\cdots +\mu _{n_{+}}\right) +\left(
n-n_{+}\right) +\left( \mu _{n_{+}+1}+\cdots +\mu _{n}\right) \\
&=&\left( e-1\right) \left( \mu _{1}+\mu _{2}+\cdots +\mu _{n_{+}}\right)
+\left( n-n_{+}\right) +\dsum\limits_{i=1}^{n}\mu _{i} \\
&=&\frac{1}{2}E_{D}\left( G\right) \left( e-1\right) +n-n_{+}.
\end{eqnarray*}

\textbf{Upper bound:}\textit{\ }Since $f\left( x\right) =e^{x}$
monotonically increases in the interval $\left( -\infty ,+\infty \right) $,
we obtain%
\begin{eqnarray*}
DEE\left( G\right) &=&\dsum\limits_{i=1}^{n}e^{\mu _{i}}\leq
n-n_{+}+\dsum\limits_{i=1}^{n_{+}}e^{\mu _{i}} \\
&=&n-n_{+}+\dsum\limits_{i=1}^{n_{+}}\dsum\limits_{k\geq 0}\frac{\left( \mu
_{i}\right) ^{k}}{k!} \\
&=&n+\dsum\limits_{k\geq 1}\frac{1}{k!}\dsum\limits_{i=1}^{n_{+}}\left( \mu
_{i}\right) ^{k}
\end{eqnarray*}%
and%
\begin{eqnarray*}
DEE\left( G\right) &\leq &n+\dsum\limits_{k\geq 1}\frac{1}{k!}\left[
\dsum\limits_{i=1}^{n_{+}}\left( \mu _{i}\right) \right] ^{k} \\
&=&n-1+e^{\frac{E_{D}\left( G\right) }{2}}.
\end{eqnarray*}

It is easy to see that the equality holds on both sides of (19) if and only
if $E_{D}\left( G\right) =0.$ Since $G$ is a connected graph this only
happens in the case of $G=K_{1}.$

\bigskip

\textbf{Remark 3.1 }It is clear that the upper bound in (19) is better than
the upper bound in (6). Moreover, the lower bounds in (7), (14) and (19) are
nicer than the lower bound (4). For instance, let $G$ be the path graph $%
P_{4}.$ Then the bounds in (7) (when $\alpha =1$ and $t=2$ ) , (14), (19)
and (4) give \linebreak $DEE\left( P_{4}\right) \geq 175.069$, $DEE\left(
P_{4}\right) \geq 92.028$, $DEE\left( P_{4}\right) \geq 11.870$ and
\linebreak $DEE\left( P_{4}\right) \geq 5.291$, respectively. Although the
best lower bound of this study is the lower bound in (7), we think that
readers like using the lower bound in (14) for practical purposes.

\bigskip

\bigskip


\bigskip

\begin{thebibliography}{99}
\bibitem{1} F. Buckley, F. Harary, \textit{Distance in Graphs, }%
(Addison-Wesley, Red-wood, 1990).

\bibitem{2} D. Cvetkovi\'{c}, M. Doob, H. Sachs, \textit{Spectra of
Graphs-Theory and }\linebreak \textit{Application }(Third ed. Johann
Ambrosius Bart Verlag, Heidelberg, Leipzig, 1995).

\bibitem{3} I. Gutman, Y. N. Yeh, S. L. Lee, Y. L. Luo, Some recent results
in the theory of the Wiener number,\textit{\ Indian J. Chem.\ }\textbf{32A}
(1993) 651-661.

\bibitem{4} E. Estrada, Characterization \ of 3D molecular structure, 
\textit{Chem. Phys. Lett. }\textbf{319} (2000) 713-718.

\bibitem{5} E. Estrada, Characterization \ of \ the folding degree of
proteins\textit{, }\linebreak \textit{Bioinformatics} \textbf{18 (}2002)
697-704.

\bibitem{6} E. Estrada, Characterization \ of \ amino acid contribution to
the folding degree of proteins\textit{, Proteins} \textbf{54 (}2004)\textbf{%
\ }727-737.

\bibitem{7} E. Estrada, J. A. Rodr\'{\i}guez-Vel\'{a}zguez, Subgraph
centrality in complex \linebreak networks, \textit{Phys. Rev. }\textbf{E 71}
(2005) 056103-056103-9.

\bibitem{8} E. Estrada, J. A. Rodr\'{\i}guez-Vel\'{a}zguez, M. Randi\'{c},
Atomic Branching in molecules, \textit{Int. J.Quantum Chem. }\textbf{106}
(2006) 823-832.

\bibitem{9} J. A. \ De la Pe\~{n}a, I. Gutman, J. Rada, Estimating the
Estrada index\textit{,} \textit{Linear Algebra Appl.} \textbf{427} (2007)
70-76.

\bibitem{10} J. P. Liu, B. L. Liu, Bounds of the Estrada index of graphs, 
\textit{Appl. Math. J. Chinese Univ.} \textbf{25 (3) }(2010) 325-330\textbf{.%
}

\bibitem{11} H. Deng, S. Radenkovi\'{c}, I. Gutman, \textit{The Estrada index%
}, in: D. Cvetkovic, I. Gutman (Eds.) \textit{Applications of Graph Spectra
\ }(Math. Inst., Belgrade, 2009) 123-140.

\bibitem{12} K. Ch. Das, S. G. Lee, On the Estrada index conjecture\textit{,}
\textit{Linear Algebra Appl. }\textbf{431} (2009) 1351-1359.

\bibitem{13} G. Indulal , I. Gutman, A. Vijaykumar, On the distance energy
of a graph\textit{,} \textit{MATCH Commun. Math. Comput. Chem.} \textbf{60}
(2008) 461-472.

\bibitem{14} G. Indulal, Sharp bounds on the distance spectral radius and
the distance energy of graphs, \textit{Linear Algebra Appl.} \textbf{430}
(2009) 106-113.

\bibitem{15} A. D. G\"{u}ng\"{o}r, \c{S}. B. Bozkurt,\textit{\ }On the
distance Estrada index of graphs\textit{, Hacettepe J. Math. Stat. }\textbf{%
38 (3)} (2009) 277-283.

\bibitem{16} A. D. G\"{u}ng\"{o}r, \c{S}. B. Bozkurt, On the distance
spectral radius and the \linebreak distance energy of graphs\textit{,} 
\textit{Linear and Multilinear Algebra}\textbf{\ 59 (4)} (2011) 365-370.

\bibitem{17} B. Zhou, N. Trinajsti\'{c}, On the largest eigenvalue of the
distance matrix of a connected graph, \textit{Chem. Phys. Lett.} \textbf{447}
(2007) 384-387.
\end{thebibliography}
\end{document}